\documentclass[12pt]{amsart}
\usepackage{amsmath,amssymb,amsfonts}
\usepackage{mathtools}
\usepackage[shortlabels]{enumitem}
\setlist{leftmargin=*}

\newcommand{\0}{\emptyset}

\makeatletter
\def\@tocline#1#2#3#4#5#6#7{\relax
  \ifnum #1>\c@tocdepth 
  \else
    \par \addpenalty\@secpenalty\addvspace{#2}%
    \@ifempty{#4}{%
      \@tempdima\csname r@tocindent\number#1\endcsname\relax
    }{%
      \@tempdima#4\relax
    }%
    \parindent\z@ \leftskip#3\relax \advance\leftskip\@tempdima\relax
    \rightskip\@pnumwidth plus4em \parfillskip-\@pnumwidth
    #5\leavevmode\hskip-\@tempdima
      \ifcase #1
       \or\or \hskip 1em \or \hskip 2em \else \hskip 3em \fi%
      #6\nobreak\relax
    \dotfill\hbox to\@pnumwidth{\@tocpagenum{#7}}\par
    \nobreak
    \endgroup
  \fi}
\makeatother

\usepackage[sorted]{amsrefs}

\newtheorem{thm}{Theorem}[section]
\newtheorem{remark}[thm]{Remark}

\newtheorem{lemma}[thm]{Lemma}
\newtheorem{prop}[thm]{Proposition}
\newtheorem{claim}[thm]{Claim}
\newtheorem{fact}[thm]{Fact}
\newtheorem*{thm*}{Theorem}
\theoremstyle{definition}
\newtheorem{defn}[thm]{Definition}
\newtheorem{rem}[thm]{Remark}

\numberwithin{equation}{section}

\newcommand{\NN}{\mathbb N}

\newcommand{\la}{\langle}
\newcommand{\ra}{\rangle}
\def\tp{{\mathrm{tp}}}
\gdef\CL{\mathcal{L}}

\def\CU{\mathcal U}

\newcommand{\sub}{\subseteq}

\usepackage[colorlinks=true, linkcolor=blue]{hyperref}

\usepackage[textsize=tiny]{todonotes}

\author{Ya'acov Peterzil} \thanks{The first author partially supported by ISF grant 290/19}
\address{University of Haifa}
\email{kobi@math.haifa.ac.il}
\author{Anand Pillay}\thanks{The second author partially supported by NSF grants DMS-1665035, DMS-1760212, and DMS-2054271}
\address{University of Notre Dame}
\email{Anand.Pillay.3@nd.edu}
\author {Fran\c{c}oise Point}
\address{University of Mons}
\email{Francoise.Point@umons.ac.be}
\begin{document}
\title[On definable groups and $D$-groups]{On definable groups and D-groups in certain fields with a generic derivation}

\begin{abstract}
We continue our study from \cite{PPP} of finite dimensional definable groups in models of the theory $T_{\partial}$, the model companion of an o-minimal $\CL$-theory $T$ expanded by a generic derivation $\partial$  as in \cite{F-K}.

 We generalize Buium's notion of an algebraic $D$-group to $\CL$-definable $D$-groups, namely $(G,s)$, where  $G$ is a $\CL$-definable group in a model of $T$, and $s:G\to \tau(G)$ is an $\CL$-definable group section.
Our main theorem says that every definable group of finite dimension in a model of $T_\partial$ is definably isomorphic to a group of the form  $$(G,s)^\partial=\{g\in G:s(g)=\nabla g\},$$ for some $\CL$-definable $D$-group $(G,s)$ (where $\nabla(g)=(g,\partial g)$).

We obtain analogous results when $T$ is either the theory of $p$-adically closed fields or the theory of pseudo-finite fields of characteristic $0$.

\end{abstract}
\maketitle
\section{Introduction}
In \cite{PPP} we initiated a study of definable groups in CODF (closed ordered differential fields, see \cite{Singer}), and more generally in  differential expansions of
o-minimal structures, $p$-adically closed fields, pseudo-finite fields of characteristic $0$, or topological fields which are models of an open theory (as in \cite{K-P}).


In all of the above settings we start with a suitable theory $T$ in a language $\CL$, where $T$  expands the theory of fields. We add a symbol $\partial$ to the language to
get $\CL_\partial=\CL\cup \{\partial\}$.  The $\CL_\partial$-theory  $T\cup$ ``$\partial$ is a (compatible) derivation'' will have a model companion which we call $T_\partial$.

 The main theorems in \cite{PPP} said that in all of these cases, if $\Gamma$ is a finite dimensional group in a model of $T_\partial$ then there is an
 $\CL$-definable group $G$ and an $\CL_{\partial}$-definable group embedding of $\Gamma$ into $G$.

Here, we mostly follow the setting suggested by Fornasiero and Kaplan, \cite{F-K}, where we start with an $\CL$-theory $T$ of an o-minimal expansion of a real closed field $K$, expand it in the language $\CL_\partial$ to the theory $T^*$ of a  $T$-compatible derivation $\partial$, and let $T_\partial$ be the model companion
of $T^*$.

In \cite{Buium}, Buium introduced the notion of an algebraic $D$-group, namely a pair $(G,s)$, where $G$ is an algebraic group and $s:G\to \tau(G)$ a rational group section
 into the prolongation of $G$. In the setting of DCF$_0$ (differentially closed fields of characteristic zero), it was shown, see  \cite[Corollary 4.2]{Pillay-foundational} and \cite{Buium},  that every finite-dimensional definable group 
 is definably isomorphic to
  $$(G,s)^\partial=\{g\in G:s(g)=\nabla(g)\},$$ ($\nabla(g)=(g,\partial g)$).

  Our goal here is to obtain analogous tools and theorems in the setting of $T_\partial$. We first associate to every $\CL$-definable $C^1$-manifold $V$, with respect to $K$, its prolongation, the bundle $\tau(V)$. We then note, as in the algebraic case,  that when $G$ is an $\CL$-definable group over a differentially closed subfield of $K$  then so is $\tau(G)$, and the projection $\pi:\tau(G)\to G$ a group homomorphism. An $\CL$-definable $D$-group is then a pair $(G,s)$ with $G$ an $\CL$-definable group and $s:G\to \tau(G)$ an $\CL$-definable group section. Our main theorem, see Theorem \ref{main D}, is:
  \begin{thm*}Let $\Gamma$ be a finite dimensional $\CL_\partial$-definable group in a model of $T_\partial$. Then there exists an $\CL$-definable $D$-group $(G,s)$ such that $\Gamma$ is definably isomorphic to
  $$(G,s)^\partial=\{g\in G:s(g)=\nabla(g)\}.$$
  \end{thm*}

When T is a model complete theory of large fields in the language of fields (plus maybe constants), Tressl, \cite{Tressl},
shows that the theory of models of T equipped with a derivation has a model companion.
He also gave a uniform (in T) axiomatization of the model companion.
 Here we treat two special cases: the case of $p$-adically closed fields and of pseudo-finite fields. We develop
  the notions of $\tau(G)$ and $(G,s)$, for a $\CL$-definable group $G$ and prove the exact analogue of the above theorem for $T_\partial$-definable groups
  (see Theorem \ref{main padic}). Along the way we prove a $p$-adic analogue of an o-minimal theorem of Fornasiero and Kaplan (see Appendix, Proposition \ref{appendix-prop}).

  When $K$ is a  pseudo-finite field we prove that every $\CL_\partial$-definable group $\Gamma$ is isogenous to $H_0\cap (H,s)^\partial$, where $(H,s)$ is an algebraic $D$-group over $K$ and
  $H_0$  a finite index subgroup of $H$ (see Theorem \ref{main bounde pac}).

  \begin{rem} The case of an arbitrary (not necessarily finite dimensional) $\CL_\partial$-definable group will be treated in a subsequent paper jointly with Silvain Rideau-Kikuchi.\end{rem}

\subsection{Preliminaries}

We refer to Section 2 of \cite{PPP} for all conventions and basic notions. Briefly, we always work in a sufficiently saturated structure and  use the fact that o-minimal structures (and later, $p$-adically closed fields) are geometric structures, to define $\dim_\CL(a/k)$ as the $acl_\CL$-dimension of $a$ over $k$. The dimension of an $\CL$-definable set $X\sub K^n$ is defined as the maximal $\dim_\CL(a/k)$, for $a\in X$ (or equivalently via cell decomposition).
If we have $\dim(a/B)=\dim X$, for $a\in X$ an $\CL$-definable set over $B$ (written also as $\CL(B)$-definable), then we say that {\em $a$ is generic in $X$ over $B$}.

For a tuple $a=(a_1,\ldots, a_n)$, we let $\partial a=(\partial a_1,\ldots, \partial a_n)$. To define the $\CL_\partial$-dimension, for $a\in K^n$ and $k\sub K$ a differential subfield,
we let $\dim_\partial(a/k)=\dim_{\CL}(a,\partial a, \ldots, \partial^n a,\ldots/k)$ (possibly infinite).
The $\CL_\partial$-dimension of a $\CL_\partial$-definable set $X\sub K^n$ over $k$ is the maximum $\dim_\partial(a/k)$, as $a$ varies in $X$.

\section{Manifolds, tangent spaces and tangent bundles}


We fix an o-minimal expansion of a real closed field $K$ in a language $\CL$. All definability in this section is in the o-minimal structure.

We first recall the basic definition of a differentiable manifold and its tangent bundle in the o-minimal setting (for differentiability in this context, see
\cite[Section 7]{vdD-book}).

\vspace{.2cm}

 \noindent{\bf Notation} Let $U\sub K^r\times K^n$ be an open definable set,  and   $f:U\to K^m$ a definable $C^1$-map, written as $f(x,y)$, $f=(f_1,\ldots, f_m)$. Given $(a,b)\in U$,
 we let $(D_x
 f)_{(a,b)}:K^r\to K^m$, and $(D_y f)_{(a,b)}:K^n\to K^m$
 denote the corresponding $K$-linear maps defined as follows:  $(D_x f)_{(a,b)}$ is the $m\times r$ matrix of partial derivatives
 $$\left (\frac{\partial f_i}{\partial x_j}(a,b)\right )_{1\leq i\leq m, 1\leq j\leq r},$$
 and $(D_y f)_{(a,b)}$ is the $m\times n$ matrix
 $$\left (\frac{\partial f_i}{\partial y_t}(a,b)\right )_{ 1\leq i\leq m, 1\leq t\leq n}.$$ Then, $(Df)_{(a,b)}$ is the $m\times (r+n)$-matrix
 $$\left ( (D_xf)_{(a,b)}, (D_y f)_{(a,b})\right ).$$

For a $C^1$ map $f:V\to W$ between open subsets of $K^n$ and $K^m$, respectively, we write $Df:V\times K^n\to W\times K^m$, for the map $$(a,u)\mapsto (f(a),(Df)_a\cdot
u),$$ where $(Df)_a\cdot u=(\sum_{i=1}^n \frac{\partial f_j}{\partial x_i} (a) u_i)_{j=1}^m.$

\vspace{.2cm}

\subsection{Definable manifolds and their tangent bundles}

\begin{defn}
{\em An $\CL$-definable $C^1$ manifold of dimension $r$},  with respect to $K$, is a topological Hausdorff space $M$, together with a finite open cover $M=\bigcup_{i=1}^n
W_i$, and homeomorphisms $\phi_i:V_i\to W_i$, where $V_i\sub K^r$ is a definable open set, such that $V_{i,j}=\phi_i^{-1}(W_i\cap W_j)$ is a definable open subset of $V_i$,
and
each map $\phi_{i.j}=\phi_j^{-1}\circ \phi_i:V_{i,j}\to V_{j,i}$ is a definable $C^1$-map (between definable open subsets of $K^r$).

The collection $\{(V_i,W_i,\phi_i)_{i\in I}\}$ is {\em an atlas} for $M$.
\vspace{.2cm}

\end{defn}

Thus, we may identify $M$ with the quotient of the disjoint union $\bigsqcup_i V_i$ by the equivalence relation
$a \sim_M b \Leftrightarrow b=\phi_{i,j}(a).$

\begin{defn}\label{def-general manifolds} For $M$ an $\CL$-definable $C^1$-manifold of dimension $r$ given as above, we let
$T(M)$ be the quotient of $\bigsqcup_{i=1}^n V_i\times K^r$ by the equivalence relation, denoted by $\sim_{T(M)}$,  given via  the  maps: $$D\phi_{i,j}:V_{i,j}\times K^r\to
V_{j,i}\times K^r\,\,;\,\,  D\phi_{i,j}(c,u)=(\phi_{i,j}(c),D(\phi_{i,j})_c\cdot u).$$
\end{defn}

We then write $$T(M)=\bigsqcup_i V_i\times K^n\diagup \sim_{T(M)},$$ and
denote
(equivalence classes of) elements in $T(M)$ by $[a,u]$, $a\in \bigsqcup_i V_i$, $u\in K^r$.

Note that if $M=U\sub K^r$ is a definable open set with the identity atlas then $T(M)=U\times K^r$.

The following are easy to verify.

\begin{lemma} Assume that $M$ and $N$ are $\CL$-definable, $C^1$-manifolds, given by atlases $(W_i,V_i,\phi_i)_{i\in I}$ and $(U_j,Z_j,\psi_j)_{j\in J}$.  If $f:M\to N$ is a $C^1$-map
(read through the charts)
then there is a well defined continuous map $Df:T(M)\to T(N)$ satisfying, whenever the elements are in the appropriate $V_i$ and $U_j$,

$$Dh([a,u])=[(f(a),D(\psi_j^{-1}\circ f\circ \phi_i)_a\cdot u)].$$
\end{lemma}

\proof In fact, the map $(a,u)\mapsto (f(a),D(\psi_j^{-1}\circ f\circ \phi_i)_{a}\cdot u)$ induces a well defined map from $T(M)$ into $\bigsqcup_{j\in J} U_j\times K^{\dim
N}$. The quotient by $\sim_{T(N)}$ gives the desired map.\qed

\begin{lemma}\label{properties of T}
\begin{enumerate}
\item For $M,N$ $\CL$-definable, $C^1$-manifolds, $T(M\times N)=T(M)\times T(N)$.

\item (Chain rule) For $f:M\to N$ and $h:N\to S$ two $\CL$-definable $C^1$-maps between $\CL$-definable $C^1$ manifolds, $D(h\circ f)=D(h)\circ D(f)$.
\end{enumerate}
\end{lemma}

Summarizing, we have:
\begin{lemma}\label{category T} $(T,D)$ is a functor from the category of definable $C^1$-manifolds to the category of definable $C^0$-manifolds. It moreover preserves
products.\end{lemma}
\section{Adding a derivation}

Let $T$ be a complete, model complete theory of an o-minimal expansion of a  real closed field $K$, in a language $\CL$. The following definition is due to Fornasiero and Kaplan,
\cite{F-K}.
\begin{defn}\label{def-compatible}
A derivation $\partial:K\to K$ is called {\em $T$-compatible} if for every $\emptyset$-definable $C^1$ map $f:U\to K^n$, for $U\sub K^m$ open, for all $a\in U$, we have
$$\partial f(a)=Df_a \cdot \partial a.$$
(Here $\partial (a_1,\ldots, a_m)=(\partial a_1,\ldots, \partial a_m)^t$).
\end{defn}

Fornasiero-Kaplan, \cite{F-K}, note that the $\CL_\partial$ theory $T\cup$ ``$\partial$ is a compatible derivation'' has a model companion, which we call $T_\partial$.
We assume from now on  that $\partial$ is a  $T$-compatible derivation on $K$, and work in models of $T_\partial$.
See \cite[Proposition 2.8, Lemma 2.9]{F-K} for instances where the compatibility condition holds.

We observe:

 \begin{claim} \label{claim0}
 Assume that $M=\bigsqcup_i V_i/\sim_M$ is a $\0$-definable manifold. Then, for $a\in M$, $\partial a$ is a well defined element of $T(M)_a$. Namely if $a_i\sim_{M} a_j$
 then $(a_i,\partial a_i) \sim_{T(M)}(a_j,\partial a_j)$.
  \end{claim}
 \proof This  is easy to verify, using the compatibility of $\partial$.\qed

\subsection{The definition of $f^{\partial}$ on an open set}

The following theorem of Fornasiero and Kaplan, which follows easily from their \cite[Lemma A.3]{F-K} plays an important role here:
In the Appendix we prove the analogous result, Proposition \ref{appendix-prop}, for $p$-adically closed fields, and the proof could be modified to give an alternative proof in the o-minimal setting as well.
\begin{fact}\label{F-K A.3} Assume that $g:W\to K^r$ is an $\CL(\0)$-definable partial function on some open $W\sub K^n\times K^m$, and $b\in \pi_2(W)\sub K^m$ is
$dcl_{\CL}$-independent.
If $g(x,b)$ is a $C^1$-map on $W^b=\{a\in K^n:(a,b)\in W\}$  then for every $a\in W^b$, the function $g$ is a $C^1$-function (of all variables) at $(a,b)$.\end{fact}

As a corollary, one obtains:\begin{fact}\label{F-K A.3.1}
If $f(x)$ is an $\CL(A)$-definable function on an open subset of $K^n$ then there is an $dcl_{\CL(\0)}$-independent tuple $b \sub A$, and a $\CL(\0)$-definable
$C^1$-function $g(x,y)$ on an open subset of $K^n\times K^{|b|}$ such that $f(x)=g(x,b)$. \end{fact}

\begin{defn} For $U\sub K^n$ open and $f:U\to K^r$ an $\CL$-definable $C^1$-map (possibly over additional parameters), let
$$f^\partial(a)=\partial f(a)-(Df)_a\partial a.$$\end{defn}

Notice that if $f$ is $\0$-definable then $f^\partial(a)=0$.
 For the following see also \cite[Lemma 2.12]{F-K}.

\begin{lemma}\label{partial f definable}If $f:U\to K^r$ is an $\CL$-definable $C^1$ map,  over a differential field $k$, then $f^\partial$ is $\CL$-definable over $k$, and
continuous on $U$.\end{lemma}
\proof  By Fact \ref{F-K A.3.1}, we may write $f(x)=g(x,b)$, for $b\in K^m$ which is $\CL(\0)$-independent, and $g$ which is a $C^1$ map,  $\CL(\0)$-definable. By the
compatibility of $\partial$,
$$\partial f(a)=\partial g(a,b)=(Dg)_{(a,b)}(\partial a,\partial b)=$$
$$=(D_xg)_{(a,b)}\partial a+(D_yg)_{(a,b)}\partial b=(Df)_a\partial a+(D_y g)_{(a,b)}\partial b.$$

It follows that $f^\partial(a)=\partial f(a)-(Df)_a\partial a=(D_yg)_{(a,b)}\partial b,$ and since $b\in k$ then so is $\partial b$. Also, because $g$ is a $C^1$-function,
$f^{\partial}$ is continuous.\qed

\begin{remark}\label{remark polynomials} When $p=\sum_m a_m x^m$ is a polynomial over $k$, then $p^\partial(x)$ is a polynomial over $k$ of the same degree:
$$p^\partial(x)=\sum_m \partial a_m x^m.$$

\end{remark}
\vspace{.2cm}

For $a\in K^n$, we let $\nabla(a)=(a,\partial a)$, and for $r\in \mathbb N$, $\nabla^r(a)=(a,\partial a,\ldots, \partial^r a).$
We also need the following:
\begin{lemma} \label{dcl and partial}Assume that $k\sub K$ is a differential field, $a\in K^m,c\in K^n$ and $c\in dcl_{\CL}(k,a)$. Then $\nabla(c)\in
dcl_{\CL}(k,\nabla(a))$.
If in addition $c$ and $a$ are $\CL$-interdefinable over $k$ then $\nabla(a)$ and $\nabla(c)$ are $\CL$-interdefinable over $k$.
\end{lemma}
\proof Assume first that $a$ is $\CL$-generic in $K^m$ over $k$. Then, $c=f(a)$ for $f$ an $\CL$-definable over $k$ and $C^1$ at $a$. We have $\partial f(a)=(Df)_a \partial
a+f^\partial(a),$ where,   by Lemma \ref{partial f definable}, $f^\partial(x)$ is $\CL$-definable over $k$. So, if we let
$h(x,u)=(f(x),(Df)_x u+f^{\partial}(x))$, then $h(\nabla (a))=\nabla(c)$, so $\nabla(c)\in dcl_{\CL}(k,a)$.

Given a general $a\in K^m$, we can write it, up to permutation of coordinates, as $(a_1,a_2)$ where $a_1\in K^{m_1}$ is $\CL$-generic over $k$ and $a_2\in dcl_{\CL}(a_1)$.
Then $c\in dcl_{\CL}(k,a_1)$, so by what we saw $\nabla(c)\in dcl_{\CL}(k,\nabla(a_1))\sub dcl_{\CL}(k,\nabla(a))$.

Finally, it clearly follows that if $a$ and $c$ are $\CL$-interdefinable over $k$ then so are $\nabla(a)$ and $\nabla(c)$.\qed







\subsection{Prolongation of functions on open sets}

Here and below, we make use of Marker's account, \cite{Marker-Manin}, of prolongations in the algebraic setting.
 \begin{defn} For $U\sub K^r$ open and $f:U\to K^n$ an $\CL$-definable $C^1$-map, we let
 $\tau(f):U\times K^r\to K^n\times K^n$ be defined as $$\tau(f)(a,u)=(f(a), (Df)_a\cdot u+f^{\partial}(a))=(f(a),(Df)_a\cdot (u-\partial a)+\partial f(a)).$$
\end{defn}

Using Lemma \ref{partial f definable} (the $\CL$-definability of $f^\partial$) and the $\CL$-definability of $Df$,
\begin{lemma} \label{tau f}If $f$ is a $C^1$-map, $\CL$-definable over a differential field $k$ then $\tau(f)$ is continuous and $\CL$-definable over the same $k$.

\end{lemma}
Using the second equality in the definition of $\tau(f)$ and the chain rule for $D$, we immediately obtain:
 \begin{lemma}\label{composition open} If $f:U\to V$ and $h:V\to W$ are definable $C^1$-functions on open sets then
 $$\tau(h\circ f)=\tau(h)\circ \tau(f).$$\end{lemma}






 \subsection{The definition of $\tau(M)$ and $\tau(f)$ for definable manifolds}

 \begin{defn} Assume that $M=\bigsqcup_i V_i/\sim_M$ is an $\CL$-definable $C^1$-manifold of dimension $n$.
 Then the prolongation of $M$ is defined as:
 $$\tau(M):=\bigsqcup_i V_i\times K^n/\sim_{\tau(M)},$$
where $(a_i,u)\sim_{\tau (M)}(a_j,v)$ if $\tau(\phi_{i,j})(a_i,u)=(a_j,v)$.
\end{defn}

By Lemma \ref{tau f}, $\tau(M)$ is an $\CL$-definable $C^0$-manifold.

The following is easy to verify.
 \begin{lemma} Assume that  $M=\bigsqcup_i V_i/\sim_M$ is an $\CL$-definable $C^1$-manifold. Then  $$  (a_i,u)\sim_{T(M)} (a_j,v) \Leftrightarrow (a_i,u+\partial
 a_i)\sim_{\tau(M)} (a_j,v+\partial a_j).$$

 In particular,
 the map $$\sigma_M:(a,u)\mapsto (a,u+\partial a)$$ induces  a well defined  $\CL_{\partial}$-definable bijection \textbf{over $M$},  between $T(M)$ and $\tau(M)$.

 \end{lemma}

Using the above lemma, we see that for $a\in M$, the element $(a,\partial a)\in \tau(M)$ is well defined (e.g., as $\sigma_M(a,0)$). We thus have a well defined map
$\nabla:M\to \tau(M)$, given in coordinates by $\nabla_M(a)=(a,\partial a)$.

\begin{defn} Assume that  $M$ and $N$ are $\CL$-definable $C^1$-manifolds, $f:M\to N$ an $\CL$-definable $C^1$ map. Then {\em the prolongation of $f$},  $\tau(f):\tau(M)\to \tau(N)$, is defined by
$$\tau(f):=\sigma_N \circ Df\circ \sigma_M^{-1}.$$
\end{defn}

The following is easy to verify:
\begin{lemma}\label{prolongation}
Assume that $M$ and $N$ as above are given via the atlases $\{(V_i,W_i,\phi_i)_{i\in I}\}$ and $\{(U_j,Z_j,\psi_j)\}$, respectively, with $\dim M=r$ and $\dim N=n$. If $f:M\to N$ is an $\CL$-definable $C^1$-map then, for $(a,u)\in V_i\times K^r$, we have
$$\tau(f)([a,u])=[\tau(\psi_j^{-1}\circ f\circ \phi_i)(a,u)].$$
\end{lemma}

\begin{lemma}\label{properties of tau}Let $M,N$ be $\CL$-definable $C^1$-manifolds.  \begin{enumerate}\item If $f:M\to N$ is $\CL$-definable over a differential field $k$ then so is $\tau(f):\tau(M)\to \tau(N)$, and
$\tau(f)$ is continuous.

\item If $f:M\to N$ and $h:N\to S$ are $\CL$-definable $C^1$ maps between $\CL$-definable $C^1$-manifolds then $\tau(h\circ f)=\tau(h)\circ \tau(f)$.

\item We have $\tau(M\times N)=\tau(M)\times \tau (N)$. Moreover, if $\pi_1:M\times N\to N$ and $\pi_1:\tau(M)\times \tau(N)\to
    \tau(M)$ are the projection maps on the first coordinates then $\tau(\pi_1)=\pi_1\circ \tau.$

\item We have $\nabla_N \circ f=\tau(f)\circ \nabla_M$.
\end{enumerate}
\end{lemma}
\proof
(1) By Lemma \ref{prolongation}, the result reduces to the $\CL$-definability of each $\tau(\psi_j^{-1}\circ f\circ \phi_i)$, and therefore follows from Lemma
\ref{partial f definable}.
(2) follows from Lemma \ref{composition open}. (3) and (4) are easy to verify.\qed

As a corollary we have:
\begin{lemma}\label{category tau}  $\tau$ is a functor from the category of definable $C^1$-manifolds to definable $C^0$ manifolds, which moreover preserves
products.\end{lemma}

\section{$\CL_\partial$-definable groups}

\subsection{Prolongation of $\CL$-definable groups, $D$-groups and Nash $D$-groups}

Let $G$ be an $\CL$-definable group, so by \cite{Pillay-groupsandfields}, it admits the structure of an $\CL$-definable $C^1$-manifold. To be precise, this is commented in \cite[Remark
2.6]{Pillay-groupsandfields}, for a structure over the reals whose functions are piecewise analytic, but the same remark holds in the $C^1$ category since definable functions in
o-minimal structure over real closed fields are piecewise $C^1$, \cite[Theorem 6.3.2]{vdD-book}.

By purely categorical reasons, using Lemma \ref{category T} and Lemma \ref{category tau} we have (see \cite[section 2]{Marker-Manin} for the same construction algebraic groups):

\begin{lemma} If $G$ is a definable group and $m:G\times G\to G$ is the group product then $$\la T(G);Dm\ra \,\,\, \mbox{ and } \la \tau(G); \tau(m)\rangle$$ are
$\CL$-definable $C^0$-groups, and the function $[a,u]\mapsto a$ is in both cases an $\CL$-definable group homomorphism from $T(G)$ and $\tau(G)$ onto $G$.

The map $a\mapsto [a,0]: G\to T(G)$ is an $\CL$-definable group section and  $\nabla_G:G\to \tau(G)$ is an $\CL_{\partial}$-definable group section.
\end{lemma}

\begin{defn} Assume that $G$ is an  $\CL$-definable group, and there exists  an $\CL$-definable group section $s:G\to \tau(G)$.
Then the pair $(G,s)$ is called {\em an $\CL$-definable $D$-group}.

\begin{rem}When $T=RCF$ is the theory of real closed fields,  every definable group admits the structure of a Nash group with respect to $K$. Namely, the underlying manifold and
group operations are semialgebraic over $K$ and either real analytic, when $K=\mathbb R$, or $C^\infty$ in general
(see discussion in \cite{HP-affineNash}, based on \cite{BCR}).
In this case every definable homomorphism between such
groups is a Nash map, thus $\pi:T(G)\to G$ and $\pi:\tau(G)\to G$ are Nash maps, and an $\CL$-definable section $s:G\to \tau(G)$ is a Nash map. We call a $D$-group $(G,s)$
in this case {\em a Nash $D$-group}.
\end{rem}
\end{defn}
Our goal is to prove:
\begin{thm}\label{main D} Let $T_\partial$ be the model companion of a complete, model complete, o-minimal theory $T$, with  a $T$-compatible derivation $\partial$. Assume that $\Gamma$ is an
$\CL_{\partial}$-definable group of finite $\CL_\partial$-dimension. Then there exists an $\CL$-definable $D$-group $(G,s)$ and an $\CL_\partial$-definable group embedding $\Gamma\to G$
whose image is $$(G,s)^\partial=\{g\in G:s(g)=\nabla_G(g)\}.$$
\end{thm}

We first recall  our result from \cite{PPP}. We shall be using the following version:.

\begin{thm} \label{paper1} If $\Gamma$ is a finite dimensional $\CL_\partial$-defined group in a model of $T_\partial$ then it can be $\CL_{\partial}$-definably embedded
in an $\CL$-definable group $G\sub K^n$ such that:

(i) Every $\CL$-generic type $p\vdash G$ is realized by some $\gamma\in \Gamma$.

(ii) There are $\CL$-definable sets $X_1,\ldots, X_r\sub G$ and $\CL$-definable functions $s_i:X_i\to K^n$ such that for each $\CL$-generic $a\in X_i$, $a\in \Gamma$ iff
$\partial a=s_i(a).$ (recall that for $a=(a_1,\ldots, a_n)$, $\partial a=(\partial a_1,\ldots, \partial a_n)$.

\end{thm}

In fact, we shall prove a more precise version of Theorem \ref{main D}:

\vspace{.2cm}

\noindent {\bf Theorem }{\em  Assume that $\Gamma$ and $G$ satisfy (i) and (ii) of Theorem \ref{paper1}. If we endow $G$ with its  $C^1$-structure, then there exists an $\CL$-definable $s:G\to
\tau(G)$,  such that $\Gamma=(G, s)^\partial$, where $(G, s)^\partial=\{g\in G: s(g)=\nabla(g)\}.$}

\vspace{.2cm}

We first prove a general fact about groups in geometric structures:

\begin{prop}\label{generic S} Let $G$ be a definable group in a geometric structure and let $S\sub G$ be a definable subset. Assume that for every generic pair $(a,b)\in
S\times S$,
we have $a\cdot b\in S$ and for every generic $a\in S$ we have $a^{-1}\in S$.

Then there is a definable $S_0\sub S$ such that $S_0\cdot S_0$ is a subgroup of $G$. Moreover, $S_0$ is a large subset of both $S$ and $S_0\cdot S_0$.

\end{prop}
\proof
We let $$S_1=\{s\in S: \mbox{ the set } \{t\in S:s\cdot t\in S\, \& \,t\cdot s\in S\} \mbox { is large in } S\}.$$

By definability of dimension in geometric structures, $S_1$ is definable. By our assumptions, $S_1$ contains all generic elements of $S$, thus, by our assumptions,  $S_0:=S_1\cap S_1^{-1}$
is also large in $S$. We claim that $S_0\cdot S_0$ is a subgroup of $G$.

We need to prove that for every $a,b,c,d\in S_0$, we have $a b c^{-1} d^{-1}\in S_0\cdot S_0$.
We fix $g\in S$ generic over $a,b,c,d$, and consider $$abc^{-1}d^{-1}=(abg)(g^{-1}c^{-1}d^{-1}).$$
By assumption on $g$, $bg\in S$ and by our choice it is in fact generic in $S$ over $a,c,d$, so in particular belongs to $S_0$.
Thus, $a(bg)\in S$, and again generic there over $c,d$, so belongs to $S_0$. Similarly, $g^{-1}c^{-1}d^{-1}\in S_0$, so $abc^{-1}d^{-1}\in S_0\cdot S_0$.
Let $H:=S_0\cdot S_0$.

To see that $S_0$ is a large subset of $H$, we fix $g\in S_0$ and $h\in S$ generic over $g$ (so $h\in S_0$). Then $g h, gh^{-1}\in S$ and generic there so in $S_0$.
It follows that $g\in S_0\cdot  S_0$  and $h\in S_0^{-1}\cdot S_0=S_0\cdot S_0$.

Hence, $S_0\sub H$ and every generic $h\in S$ over $g$ is in $H$, so $S_0$ is large in $H$.\qed

We are now ready to prove Theorem \ref{main D}.

  We first apply Theorem \ref{paper1}, and deduce the existence of pairwise disjoint $\CL$-definable  $X_1,\ldots, X_r\sub G\sub K^n$, each  of dimension equal to $m=\dim
  G$, such that $X=\bigsqcup_i X_i$ is a large subset of $G$ and on each $i$, we have an $\CL$-definable $s_i:X_i\to K^n$, such that for  $g$ generic in $X_i$, we have $g\in
  \Gamma \Leftrightarrow \partial g=s_i(g)$. We let $s:X \to K^n$ be the union of the $s_i$'s.

Exactly as in \cite[Proposition 2.5]{Pillay-groupsandfields}, there exists a large $\CL$-definable set $W\sub G$, and an $\CL$-definable homeomorphism  $\sigma:V\to W$, for $V\sub K^m$
definable open,
and  $g_1,\ldots, g_k\in G$,
such that

(i) $G=\bigcup_j g_j W$.

(ii) The maps $\phi_i: V\to g_iW: x\mapsto g_i\sigma(x)$ endow $G$ with a definable $C^1$-manifold structure,
and make $G$ into a $C^1$-group.

(The original theorem is for continuous functions and a topological group, but as was remarked earlier the same proof works for $C^1$-functions and a $C^1$-group).

By intersecting $W$ with the relative interior of $X$ in $G$, we may assume that $X=W$.

\begin{claim} \label{claim hat s}There exists an $\CL$-definable $\hat s:V\to K^m$, such that for every $\CL$-generic $a\in V$, $\hat s(a)=\partial a \Leftrightarrow
s(\sigma(a))=\partial \sigma(a)$.\end{claim}
\proof Every $a\in V$ is $\CL$-interdefinable with $\sigma(a)$, so by Lemma \ref{dcl and partial}, $\nabla(a)$ and $\nabla(\sigma(a))$ are $\CL$-interdefinable over $k$.
By compactness, there exists an $\CL$-definable (partial) bijection $h:W\times K^n\to V\times K^m$, such that for each generic $a\in V$, $h(\nabla(\sigma
a))=\nabla(a)$. Let
$$\hat s(a)=\pi_2(h(\sigma(a)),s(\sigma(a))),$$where $\pi_2:V\times K^n\to K^n$ is the projection onto the second coordinate.

Now, if $s(\sigma(a))=\partial(a)$, then $(\sigma(a),s(\sigma(a)))=\nabla(\sigma(a))$, so $$\hat s (a)=\pi_2( h(\nabla \sigma(a)))=\pi_2(\nabla(a))=\partial a.$$
The converse follows from the invertibility of $h$.\qed

Going back to $G$, we now endow $G$ with a finite $C^1$-atlas $(V_i,g_iW,\phi_i)_{i\in I}$, where $V_i=V$ for all $i$,  and identify $G$ with $\bigsqcup V_i/\sim_M$.
We also identify $\Gamma$ with the group $\bigsqcup \phi_i^{-1}(\Gamma\cap g_iW)/\sim_M$.
Notice that each $g_iW/\sim_M$ is large in $G$, and  by Claim
\ref{claim hat s}, there is an $\CL$-definable $\hat s:V\to K^m$ such that, for generic $a\in V$,  $\hat s(a) =\partial a$ if and only if $s(\sigma(a))=\partial \sigma(a)$.
Thus, by our assumption, for every generic $g\in G$, $g\in \Gamma \Leftrightarrow \hat s(g)=\partial g$. For simplicity, from now on we use $s$ instead of $\hat s$ and let
$X=dom (s)$, an $\CL$-definable large subset of $G$.

Consider the $\CL$-definable $C^1$-group $\tau (G)$ as before, and the associated $\CL$-definable homomorphism $\pi:\tau(G)\to G$, together with an $\CL_\partial$-definable
group section $\nabla_G:G\to \tau(G)$.
The map $s$ can be replaced by $x\mapsto (x,s(x))$, so we may think of it as a function from $X$ into  $\tau(X)=X\times K^m$ with $\pi\circ s(x)=x$.

In addition, we still have for every generic $g\in X$,  $g\in \Gamma \Leftrightarrow s(g)=\nabla_G(g)$.
By our assumptions, every generic $\CL$-type of $X$ contains an element of $\Gamma$, hence the $\CL$-definable set  $X_0=\{x\in X:s(x)\in \tau(X)\}$ is large in $X$, so
without loss of generality, $X=X_0$.  Let $S$ be the graph of $s|X_0$.

We claim that $S$ satisfies the assumptions of Proposition \ref{generic S}: Indeed, assume that  $(a,b)$ is generic in $S^2$.
Namely, $a=(g,s(g))$ and $b=(h,s(h))$,
for $(g,h)$ generic in $X\times X$.  We need to prove that $a b\in S$.

By \cite[Lemma 6.7]{PPP}, applied to the function $(s,s):X\times X\to \tau(X\times X)$, there exists $(x,y)\in X\times X$, realizing the same $\CL$-type as
$(g,h)$ such that $\nabla_{G\times G}(x,y)=(s(x),s(y))$. But then, by our assumptions, $(x,y)\in \Gamma\times \Gamma$, so $xy\in \Gamma$. Because $xy$ is still $\CL$-generic
in $G$, we have $xy\in X$. Thus we have $$s(xy)=\nabla_G(xy)=\nabla_G(x) \nabla_G(y)=s(x) s(y)$$ (where the middle equality follows from the fact that $\nabla_G$
is a group homomorphism). Since $tp_{\CL}(x,y)=tp_{\CL}(g,h)$, we also have $s(gh)=s(g)s(h)$, hence
 $a b=(gh,s(gh))$,  is in $S$.

We similarly prove that for $a$ generic in $S$, we have $a^{-1}\in S$, thus $S$ satisfies indeed the assumption of Proposition \ref{generic S}.

Hence, there exists an $\CL$-definable $S_0\sub S$, such that $S_0$ is a large subset of the group $H=S_0\cdot S_0$. Since $S_0$ is large in $H$, for every generic
$(g,s(g))\in H$, we have $\pi^{-1}(g)\cap H$ is a singleton, which implies that $ker(\pi|H)=\{1\}$, hence  $H$ is the graph of a function. Also, since the group $\pi(H)$ is large in
$G$, it necessarily equals to $G$.

We therefore found an $\CL$-definable group-section $\hat s:G\to \tau(G)$, making $(G,\hat s)$ into a $D$-group. In addition, $x\in \Gamma \Leftrightarrow \hat
s(x)=\nabla_G(x)$, for all $x$ generic in $G$.

It is left to see that
$$\Gamma=(G,\hat s)^\partial=\{x\in G:\hat s(x)=\nabla_G(x)\}.$$

Let $X_0=\pi(S_0)$ and $\Gamma_0=X_0\cap \Gamma$. By the definition of $S$, $\Gamma=\{x\in \pi(S):S(x)=\nabla_G(x)\}$.\, so $\Gamma_0=\{x\in X_0:\hat s(x)=\nabla_G(x)\}.$
We claim that $\Gamma_0\cdot \Gamma_0=\Gamma$.

Indeed, let $\gamma\in \Gamma$, and pick $g$ generic in $X_0$ over $\gamma$. By the geometric axioms, there exists $\gamma_1\equiv_{\CL(\gamma)} g$ such that $\hat
s(\gamma_1)=\nabla_G(\gamma_1)$, namely $\gamma_1\in \Gamma_0$. It follows that $\gamma\cdot \gamma_1^{-1}$ is $\CL$-generic in $G$ over $\gamma$ and hence in $X_0$, namely
in $\Gamma_0$. Hence, $\gamma\in \Gamma_0\cdot \Gamma_0$.

It follows that for all $\gamma\in \Gamma$, we have $\hat s(\gamma)=\nabla_G(\gamma)$.  To see the converse, assume that $\hat s(x)=\nabla_G(x)$, and choose $\gamma\in
\Gamma_0$ generic over $x$. We then have $\hat s(\gamma)=\nabla_G(\gamma)$, and $x\cdot \gamma$ generic in $X_0$.
Because $\hat s$ is a homomorphism,
$$\hat s(x\gamma)=\hat s(x)\hat s(\gamma)=\nabla_G(x)\nabla_G(\gamma)=\nabla_G(x\gamma).$$ It follows that $x\gamma\in \Gamma$ and hence so is $x$. This ends the proof of
Theorem \ref{main D}.\qed

\subsection{The case of $p$-adically closed fields}

Let $K$ be a $p$-adically closed field, namely a field which is elementarily equivalent to a finite extension of $\mathbb Q_p$.
The field admits a definable valuation, which we may add to the field language and call this language $\CL$.

We shall use multiplicative notation for the valuation map $|\,|:K\to \{0\}\cup \Gamma$. Namely, $$|0|=0<\Gamma\,\, ,  |\,|:K^*\to (\Gamma,\cdot)\, \mbox{a
group homomorphism}$$ and $$\forall x,y\in K \,\,|x+y|\leq \max \{|x|,|y|\}.$$
For $a=(a_1,\ldots,a_n)\in K^n$, we write $\| h\|=\max\{|a_i|:i=1,\ldots,n\}.$ Since $K$ is a geometric structure we use  the $acl$-dimension below.

\begin{defn} For $U\sub K^m$ open, a map $f:U\to K^n$ is called {\em differentiable at $a\in U$} if there exists a $K$-linear map $T:K^m\to K^n$
such that for all $\epsilon\in \Gamma$ there is $\delta\in \Gamma$, such that for all $h\in K^m$, if $\|h\|<\delta$ then
$$\|f(a+h)-f(a)+T(h)\|<\epsilon\|h\|.$$

The linear map $T$ can be identified with $Df_a$ the $n\times m$ matrix of partial derivatives of $f$. We identify $M_{n\times m}(K)$ with $K^{n\cdot m}$.
The function $f$ is called {\em continuously differentiable on $U$}, or $C^1$, if it is differentiable on $U$ and the map $x\mapsto Df_a$ is continuous.
\end{defn}

Differentiable maps satisfy the chain rule, by the usual proof (see for example \cite[Remark 4.1]{Schneider} for a proof in $\mathbb Q_p$).

Towards our main result,
we first note that every $\CL$-definable group in $K$ can be endowed with the structure of a $C^1$-Lie group  with respect to the valued field $K$, with finitely many
charts:

 Indeed,
in \cite[Lemma 3.8]{Pillay-Qp-fields} a similar statement was proved for definable groups in $\mathbb Q_p$, based on the work \cite{Pillay-groupsandfields} from the
o-minimal setting. The following observations are needed in order to make adjustments to the $p$-adically closed setting (see \cite[Fact 3.7]{Pillay-Qp-fields}):

(i) Every definable $X\sub K^n$ can be partitioned into a finite number of definable sets $X_i$, each of which is homeomorphic by projection along certain coordinate axes
to an open subset of $K^m$ for some $m$. (see \cite[Theorem 1.1 and Section 5]{Scowcroft-vdd} for finite extensions of $\mathbb Q_p$).

\vspace{.2cm}

(ii) If $X\sub K^n$ is a definable open set and $f:X\to K$ a definable function then the set $Y$ of $x\in K^n$ such that $f$ is $C^1$ at $x$ is large in $X$,
namely $\dim(X\setminus Y)<\dim(X)$. (see \cite[Theorem 1.1']{Scowcroft-vdd}).

\vspace{.1cm}

Because $p$-adically closed fields are geometric structures, dimension is definable in parameters. In addition, so is  differentiability. Thus (i) and (ii) follow from the fact that the
statements are  true in finite extensions of $\mathbb Q_p$.

\vspace{.2cm}

Using (1) and (2) above, it indeed follows that every $\CL$-definable group in $K$ can be endowed with the structure of a $C^1$-group with respect to $K$.
(The same proof shows that  every definable group in $\mathbb Q_p$ admits a definable $p$-adic Nash group structure over $K$).

\vspace{.2cm}

We now endow $K$ with a derivation, denoted by $\partial$.
By Tressl's work, see \cite[Theorem 7.2]{Tressl}, the theory of $p$-adically closed fields with a derivation has a model companion $T_{\partial}$. In our one
derivation case (Tressl deals with several commuting derivations), one can axiomatize $T_{\partial}$ with the following geometric axioms (see for instance \cite[Fact 5.7 (ii)]{PPP}):
whenever $(V,s)$ is an irreducible  $D$-variety over $K$ with a smooth $K$-point and $U$ is a Zariski open subset of $V$ defined over $K$, then there is $a\in U(K)$ such that
$(a,s(a))=\nabla(a)$. (Recall that a $D$-variety $(V,s)$ defined over $K$ is a $K$-variety $V$ equipped with a rational section $s$ defined over $K$ from $V$ to $T(V)$ \cite[Definition
2.4]{PPP}).

Now, exactly as in the work of Fornasiero and Kaplan for real closed fields,  \cite[Lemma 2.4, Lemma 2.7, Proposition 2.8]{F-K}, the associated derivation is compatible with the theory
of $p$-adically closed fields, namely compatible with every $\CL(\0)$-definable $C^1$ map, as in  Definition
\ref{def-compatible}.
In order to develop the rest of the theory as in the o-minimal case, we prove in the Appendix (see Proposition \ref{appendix-prop}) that definable functions in $p$-adically closed fields satisfy
the analogue of Fact \ref{F-K A.3.1}:

\begin{prop} \label{appendix-prop1} Given an $\CL(\0)$-definable $W\sub K^n\times K^m$ and an $\CL(\0)$-definable
$g:W\to K$, if $(a,b)\in W$, $\dim(b/\0)=m$, $a\in Int(W^b)$ and $g(x,b)$ is a $C^1$-function on $W^b$ then $(a,b)\in Int(W)$ and $g$ is a $C^1$-function at $(a,b)$.

\end{prop}

Now, the category of $K$-differentiable manifold\textcolor{blue}{s} $M$ and their associated functors $T$ and $\tau$ can be developed identically to Sections 1 and 2.  This allows us to associate to every
definable group $G$ the definable groups $T(G)$ and $\tau(G)$, such that the natural projections onto $G$ are group
homomorphisms. If $G$ is a $C^1$-group then  $T(G)$ and $\tau(G)$ are $C^0$-groups.

By a $p$-adic $D$-group, we mean a pair $(G,s)$ where $G$ is an $\CL$-definable $C^1$-group and $s:G\to \tau(G)$ an $\CL$-definable homomorphic section (i.e., $\pi\circ s=id$).

As before, we define in models of $T_{\partial}$, given a $D$-group $(G,s)$,
$$(G,s)^\partial=\{g\in G:s(g)=\nabla_G(g)\}.$$


Using the analogue of Theorem \ref{paper1}, proven in \cite{PPP}, we may repeat the exact same proof as
in the o-minimal case to conclude:
\begin{thm}\label{main padic} Let $T$ be the theory of $p$-adically closed fields and let  $\Gamma$ be a finite dimensional  $L_{\partial}$-definable group in $K\models
T_{\partial}$. Then there exists an $\CL$-definable $D$-group $(G,s)$ such that $\Gamma$ is definably isomorphic to $(G,s)^\partial$.\end{thm}

\subsection{The case of pseudofinite fields}
Let $\CL$ be the language of rings and let $C=(c_{i,n})_{n\in \mathbb N, i<n}$ be an infinite countable set of new constants. Let $T$ be the $\CL(C)$-theory of pseudo-finite fields of characteristic $0$, namely the
theory of pseudo-algebraically closed fields plus the scheme of axioms saying, for every $n\in \NN$, that there is a unique extension of degree $n$,
and that the polynomial
$$X^n+c_{n-1,n}X^{n-1}+\cdots c_{0,n}$$ is irreducible.

  Since $T$ is a model-complete theory of large fields,
one can apply the Tressl machinery and so the theory of differential expansions of models of $T$ has a model-companion \cite[Corollary 8.4]{Tressl}, which has been axiomatized
\cite[Theorem 7.2]{Tressl} (in case of expansions by a single derivation, one obtains a geometric axiomatisation  \cite[Lemma 1.6]{BCPP}). Recall that since $T$ has almost q.e. (see
\cite[Remark 1.4 (2)]{BCPP}), the theory $T_\partial$ does too \cite[Definition 1.5, Lemma 2.3]{BCPP}, \cite[Theorem 7.2(iii)]{Tressl}.
\medskip

Let ${\mathcal U}$ be our sufficiently saturated model of $T_\partial$, a differential extension of a pseudo-finite field and let $\bar{\mathcal U}\supseteq \mathcal U$ be a saturated model of DCF$_0$ extending it.  We work over a small submodel $(K,\partial)\models T_\partial$.


We briefly review the construction of the algebraic prolongation $\tau(V)\sub \bar{\CU}^n\times \bar{\CU}^n$ of an
irreducible algebraic variety (see \cite{Marker-Manin} for details):

Assume that the ideal $I(V)$ is generated by polynomials $p_1,\ldots, p_m$, over $K$, and let $P:\bar{\CU}^n\to \bar{\CU}^m$ be the corresponding polynomial map $P(x)=(p_1(x),\ldots,
p_m(x))$.
The definition of $DP$ and $\tau(P)$ is defined as before using the formal derivative of polynomials (see also Remark \ref{remark polynomials}).
Then

$$T(V)=\{(x,u)\in \bar{\CU}^{2n}: P(x)=0\,\& \,(DP)_x\cdot u=0\},$$ and
$$\tau(V)=\{(a,u)\in \bar{\CU}^{2n}: a\in V\, \&\, \tau(P)(a,u)=0\}.$$

Both are  algebraic varieties over $K$.
For $a\in V(\bar{\CU})$,  $a\mapsto \partial a$ is a section of $\pi:\tau(V)\to V$, and we have
$$\tau(V)=\{(a,u)\in \bar{\CU}^{2n}: a\in V\, \&\, u-\partial a\in T(V)_a\}.$$ So, $\tau(V)_a$ is an affine translate of the vector space $T(V)_a\sub \bar {\CU}^n$. In particular,
$\dim(\tau(V)_a)=\dim V$.

As described in \cite{Marker-Manin}, the above constructions of $T(V)$ and $\tau(V)$ can be extended to abstract, not necessarily affine, algebraic varieties (which are
covered by finitely many affine algebraic varieties). Furthermore, if $H$ is an algebraic group, then $T(H)$ and $\tau(H)$ are algebraic groups with the property that the map
$\nabla_{H}$ is now a group morphism \cite[section 2]{Marker-Manin}.


Our goal is to prove:

\begin{thm}\label{main bounde pac} Let $T$ be the theory of pseudo-finite fields and let $\Gamma$ be a finite dimensional definable group in $K\models T_{\partial}$. Then
there exists a $K$-algebraic $D$-group $(H, s)$
such that $\Gamma$ is virtually definably isogenous over $K$ to the $K$-points of $(H,s)^\partial$.\end{thm}

By ``$\Gamma$ and $(H,s)^\partial(K)$ are virtually isogenous'' we mean the following: There exists an $\CL_{\partial}$-definable subgroup $\Gamma_0\sub \Gamma$ of finite
index and a $\CL_{\partial}$-definable homomorphism $\sigma:\Gamma_0\to H$ with finite kernel, whose image has  finite index in the $K$-points of $(H,s)^\partial$.

We first need the following lemma:

\begin{lemma}\label{isomorphism}  Let $W_1, W_2$ be irreducible algebraic  varieties over a differential field $K$, and $a, b$ generic tuples in $W_1$ and $W_2$,
respectively, over $K$.  Suppose that $a$ and $b$ are field-theoreticallly interalgebraic over $K$ (in particular, $\dim W_1=\dim W_2$), and let $W\sub W_1\times W_2$ be the
(irreducible) variety over $K$ with generic point $(a,b)$.
Then $\tau(W)_{(a,b)}$ is the graph of a bijection over $K(a,b)$,  $\alpha:\tau(W_{1})_{a}\to \tau(W_{2})_{b}$.
\end{lemma}
\proof
By dimension considerations, the affine space $\tau(W)_{(a,b)}$ projects onto both $\tau(W_{1})_{a}$ and $\tau(W_{2})_{b}$.

Fix some coordinate $a_{1}$ of the tuple $a$.  By our interalgebraicity assumption, $a_{1}$ is in the field-theoretic algebraic closure of $K(b)$.
Let $q(x)$ be the minimal monic polynomial of $a_{1}$ over $K(b)$.
$q(x) =  x^{n} + f_{n-1}(b)x^{n-1} + ... f_{1}(b)x + f_{0}(b)$, where the $f_{i}(b)$ are $K$-rational functions of the tuple $b$.
After getting rid of denominators  we can rewrite $q(x)$ as
$q_n(b)x^{n} + q_{n-1}(b)x^{n-1} + ... + q_{1}(b) x + q_{0}(b)$  where the $q_{i}$ are polynomials over $K$.

Hence  $r(x,y)$:   $q_n(y)x^{n} + q_{n-1}(y)x^{n-1} + ... + q_{0}(y)$  is a polynomial (over $K$) in $I_{K}(W)$.

Hence  $(\partial r/\partial x)(a_{1}, b)(u_{1}) +  \sum_{j}(\partial r/\partial y_{j})(a_{1},b)(v_{j}) + r^{\partial}(a_{1},b) = 0$ for
$(u_{1},...,v_{1},..,v_{j}...)\in \tau(W)_{a,b}$.

By the minimality of $r(x,b)$, we have $\partial r/\partial x(a_1,b)\neq 0$, hence
$$u_{1} =  (\sum_{j}(\partial r/\partial y_{j})(a_{1},b)(v_{j}) + r^{\partial}(a_{1},b))/ (\partial r/\partial x)(a_{1}, b)).$$ Thus, $u_1\in dcl(K,a,b,v)$. We similarly prove that $u$
and $v$ are inter-definable over $K(a,b)$. Since
$\tau(W)_{(a,b)}$ is a translate of a linear space, containing $(u,v)$ whose dimension equals $\dim W_1=\dim W_2$, it must be the graph of a bijection.\qed

We now return to the proof of Theorem \ref{main bounde pac}:

By \cite[Theorem 5.11]{PPP} and the associated Remark 4.8 there, (formulated for general large geometric fields), we may assume that $\Gamma \sub G\sub \CU^n$ for some
$\CL$-definable group $G$ such that every generic $\CL$-type of $G$ is realized by an element in $\Gamma$. Furthermore, there is a covering of $G$ by finitely many
$\CL(K)$-definable sets $X_i$, $i=1,\ldots, m$, and for each $X_i$ there is a $K$-rational function $s_i:X_i\to \CU^n$ such that for every $a\in X_i(\CU)$ which is $\CL$-generic
in $X_i$ over $K$, we have
$$a\in \Gamma \Leftrightarrow \partial(a)=s_i(a).$$
We may take each $X_i$ to be Zariski  dense in a $K$-variety $V_i$.
We are only interested in those $V_i$ whose Zariski dimension is maximal, call it $d$, so
in particular, every algebraic type in $V_i$ over $K$, of dimension $d$, is realized in $X_i$ in $\CU$, so by the axioms also realized by some $a\in X_i(\CU)$ with $\partial a=s_i(a)$, and hence, by the above, also realized in $\Gamma$.
Finally, each $X_i$ can be taken to be  the $\CU$-points of $W_i=Reg(V_i):=V_i\setminus Sing(V_i)$, namely the $\CU$-points of a smooth quasi affine $K$-variety.

We now apply \cite[Theorem C]{Hrushovski-Pillay-groupslocalfields} in the structure $\CU$: There exist a connected algebraic group $H$ over $K$,   $\CL$-definable subgroups
of finite index, $G_0\sub G$ and  $H_0\sub H(\CU)$, and an  $\CL$-definable surjective homomorphism $f:G_0\to H_0$ whose kernel is finite, all defined over $K$. The group
$H$, as an algebraic group over $K$, has an associated $K$-algebraic group $\tau(H)$. Notice $acl_{\CL}$ equals the field $acl$ and we have $\dim H=d$.

Let $\Gamma_0=\Gamma\cap G_0$, a subgroup of $\Gamma$ of finite index.
Since $f(G_0)$ is Zariski dense in $H$ so is $f(\Gamma_0)$. As we shall now see, we can endow $H$ with the structure of a $D$-group, $(H,s)$ such that
$$f(\Gamma_0)=\{h\in H(\CU):s(h)=\nabla_H (h)\}.$$

\begin{claim} For every $g\in \Gamma_0$, $tr.deg(\nabla_H f(g)/K)\leq \dim(H).$\end{claim}
\proof We first prove the result for $g\in \Gamma_0$ such that $tr.deg(g/K)=d$.

Since $tr.deg(g/K)=d$, there exists a $K$-algebraic quasi-affine variety $W_i$ as above such that $g$ is generic in $W_i$ over $K$. The $\CL$-definable function $f$
takes values in $H$, and because $acl_{\CL}$ is the same as the field $acl$, there exists an algebraic correspondence $C_i\sub W_i\times H$ over $K$, such that $(g,f(g))$ is
field-generic in $C_i$. It follows that $(\nabla_{W_i}(g),\nabla_H f(g))\in \tau(C_i)\sub \tau(W_i)\times \tau(H)$.

By Lemma \ref{isomorphism},  $\tau(C_i)_{(g,f(g))}$ induces an (algebraic) bijection over $K$ between $\tau(W_i)_g$ and $\tau(H)_{f(g)}$.
In particular, $\nabla_{W_i}(g)$ and $\nabla_H f(g)$ are interalgebraic over $K$. By our construction, $\nabla_{W_i}(g)=s_i(g)$, hence $g$ and $\nabla_H f(g)$ are
interalgebraic over $K$ (notice that $g$ is algebraic over $\nabla_{W_i}(g)$). Hence, $tr.deg(\nabla_H f(g)/K)=tr.deg(g/K)=d$.

Assume now that $g$ is an arbitrary element of $\Gamma_0$, and let $h\in \Gamma_0$ be such that
$tr.deg(g/K,\nabla_Hf(g))=d$.

Since $tr.deg(hg/K)=d$ and $hg\in \Gamma_0$,  it follows from the above that $$d=tr.deg(\nabla_{H}f(hg)/K))=tr.deg(\nabla_H f(h)\cdot \nabla_H f(g)/K),$$
and therefore $$tr.deg(\nabla_Hf(h) \cdot \nabla_H f(g)/\nabla_H f(h), K)\leq d.$$ The elements $\nabla_Hf(h)\cdot \nabla_H f(g)$ and $\nabla_H f(g)$ are interalgebraic over
$K$ and $\nabla_Hf(h)$, thus $tr.deg(\nabla_H f(g)/\nabla_H f(h),K)\leq d$.

We know that $h$ and $\nabla_{X_i}(h)$ (and hence also $\nabla_H f(h)$) are interalegbraic over $K$ (as witnessed by $s_i$), and as $h$ and $\nabla_H(f(g))$ are independent
over $K$, it follows that $\nabla_H f(h)$ and $\nabla_H f(g)$ are independent over $K$. Therefore,
$$tr.deg(\nabla_H f(g)/K)=tr.deg(\nabla_Hf(g)/\nabla_H f(h),K)\leq d.$$ \qed

We now consider the subgroup  $\nabla_H f(\Gamma_0)$ of $\tau(H)$ and let $S\sub \tau(H)$ be its Zariski closure, an algebraic subgroup of $\tau(H)$.
By the claim above, $\dim(S)\leq \dim H$, but since $S$ contains $\nabla_H f(h)$ for $\CL$-generic $h\in \Gamma_0$, we have $\dim(S)=\dim(H)$. Consider the projection
$\pi:\tau(H)\to H$, a group homomorphism, and its restriction to $S$. Since $H$ is connected, we have $\pi(S)=H$, and hence $ker(\pi)\cap S$ is a finite subgroup of $\tau(H)_e$. However,
$\tau(H)_e=T(H)_e$ is a vector space over $K$, a field of characteristic $0$, thus torsion-free. Hence, $ker(\pi)\cap S$ is trivial so $\pi:S\to H$ is a group isomorphism.
It follows that $S$ can be viewed as a group section $s:H\to \tau(H)$.
Since $S$ is the Zariski closure of $\nabla_H(f(\Gamma_0))$, we have for every $g\in \Gamma_0$, $\nabla_H(f(g))=s(f(g))$.

Recall that $H_0=f(G_0)$ is an $\CL$-definable subgroup of of finite index of $H(\CU)$.

\begin{claim} $f(\Gamma_0)=(H_0,s)^\partial=\{h\in H_0:\nabla_H(h)=s(h)\}.$ \end{claim}

\proof We only need to prove the $\supseteq$ inclusion.


We first prove that for every $h\in (H_0,s)^\partial$, if $tr.deg(h/K)=d$ then $h\in f(\Gamma_0)$. Indeed, since $h\in H_0$ is $\CL$-generic in $H$ over $K$, there exists $g\in
G_0$,  necessarily $\CL$-generic in $G$ over $K$, such that $h=f(g)$.  By our assumptions, there exists $g'\in \Gamma$,  such that $g'$ and $g$ realize the same $\CL$-type over
$K$.
Thus, $g'$ is in  $G_0$ so also in $\Gamma_0$. In addition, $f(g')$ and $f(g)=h$ must realize the same $\CL$-type over $K$. Because $S$ is the Zariski closure of $\nabla_H
f(\Gamma_0)$, it follows that $f(g')\in (H_0,s)^\partial$.

Since $\nabla_H(h)=s(h)$ and  $\nabla_H(f(g'))=s(f(g'))$, it follows that for every $n\in \mathbb N$, there is an $\CL(K)$-definable function $s_n$ such that
$$\nabla^n_H(h)=s_n(h)\, ,\,\nabla^n_H(f(g'))=s_n(f(g'))$$
(recall that $\nabla_H^n(g)=(g,\partial g,\ldots,\partial^n g)$). Because  $h$ and $f(g')$ realize the same $\CL$-type over $K$, we may conclude that for every $n$,
$\nabla_H^n(h)$ and $\nabla_H^n(f(g')))$ realize the same $\CL$-type over $K$, and therefore
$$\tp_{\CL(K)}(\nabla_H^{\infty}(h))=\tp_{\CL(K)} (\nabla_H^{\infty}(f(g'))).$$
($\nabla_H^\infty(g)=(g,\partial g,\ldots, \partial^n g,\ldots).$)

By \cite[7.2(iii)]{Tressl}, every $\CL_{\partial}(K)$-formula is equivalent in $T_{\partial}$ to a boolean combination of formulas of the form $(x,\partial x,\ldots,\partial^n x)\in Y$, where $Y$ is an $\CL(K)$ definable set. Thus, it follows from the above that $h$ and  $f(g')$ realize  the same $\CL_{\partial}$-type over $K$, and therefore $h\in f(\Gamma_0)$, as needed.


This proves that every  $h\in (H_0,s)^\partial$ with $tr.deg(h/K)=d$ belongs to $f(\Gamma_0)$. However, every $h\in (H_0,s)^\partial$ can be written as $h=h_1h_2$, with $h_1,h_2\in (H_0,s)^\partial$ and $tr.deg(h_1/K)=tr.deg(h_2/K)=d$. Indeed, pick $h_1\in (H_0,s)^\partial$ with $tr.deg(h_1/hK)=d$ and $h_2=h_1^{-1}h$). Thus,  every $h\in (H_0,s)^\partial$ belongs to $f(\Gamma_0)$.
This ends the proof of Theorem \ref{main bounde pac}.\qed

\vspace{.3cm}

\section{Appendix}

We fix $K$ a $p$-adically closed field. All definability below is in the language $\CL$ of $K$.
Our goal is to prove the $p$-adic analogue of Fornasiero-Kaplan's theorem (see \cite[A.3]{F-K}).
Before stating it we make some preliminary observations.





We can now state the result that we plan to prove here:

\begin{prop}\label{appendix-prop} Let $K$ be a $p$-adically closed field.
Assume that $g:W\to K^r$ is an $\CL(A)$-definable partial function on some definable $W\sub K^n\times K^m$, and $b\in \pi_2(W)\sub K^m$ is $acl_{\CL}$-independent over $A$.
If $g(x,b)$ is a $C^1$-map on the open set $W^b=\{a\in K^n:(a,b)\in W\}$  then for every $a\in W^b$, the function $g$ is a $C^1$-map (of all variables) in a neighborhood of
$(a,b)$.\end{prop}

We shall use the following three important properties of $p$-adically closed fields (as well as o-minimal structures and some other geometric structures).
\begin{fact}\label{fact-app}
Fix $A\sub K^{eq}$.
\begin{enumerate}
\item
Given $a\in K^m$ and $b\in K^n$, if $U\ni a$ is a (definable) open set in $K^m$ then there exists a definable open $V$, $a\in V\sub U$, such that $\dim(b/A[V])=\dim(b/A)$ (we use
$[V]$ for the canonical parameter of $U$), see \cite[Corollary 3.13]{HHP}, or \cite[Lemma 4.30]{Johnson}.

\item Assume that $X\sub K^{m+n}$ is definable over $A$, $a\in K^m, b\in K^n$ and $(a,b)\in X$. Assume further that $X^b=\{x\in K^m:(x,b)\in X\}$ is finite. Then there exists a
    definable open $W\ni (a,b)$ (possibly over additional parameters) such that $X\cap W$ is the graph of a definable map from $K^m$ to $K^n$ (this follows from cell decomposition in
    $p$-adically closed fields).

\item If $U\sub K^m$ is open and $f:U\to K^n$ is an $A$-definable map then $f$ is $C^1$ at every $a\in U$ with $\dim(a/A)=m$ (see \cite{Scowcroft-vdd}).
\end{enumerate}

\end{fact}

An immediate corollary of the first two is:
\begin{fact}\label{appendix-cor}
For $a\in K^m$ and $b\in K^n$ and $A\sub K^{eq}$, if $a\in acl(b,A)$ then there exists $A_1\supseteq A$ such that $a\in dcl(bA_1)$ and $\dim(b/A_1)=\dim(b/A)$.
\end{fact}

\vspace{.2cm}

\noindent {\bf Proof of Proposition \ref{appendix-prop} :}  We first prove a continuous version.

\begin{lemma}\label{appendix-lemma1} Assume that $g:W\to K^r$ is an $\CL(A)$-definable partial function on some definable $W\sub K^m\times K^n$, and $b\in \pi_2(W)\sub K^n$ is
$acl_{\CL}$-independent over $A$.
If $g(x,b)$ is continuous on the open set $W^b=\{a\in K^n:(a,b)\in W\}$  then $(a,b)\in Int(W)$ and for every $a\in X^b$, the function $g$ is continuous at $(a,b)$.
\end{lemma}
\proof We need the following claim:

\begin{claim} \label{appendix-claim} Assume that $X\sub K^m\times K^n$ is an  $A$-definable set, $(a,b)\in X$ and $b$ is generic in $K^n$ over $A$.
If $a\in Int(X^b)$ then $(a,b)\in Int(X)$.
\end{claim}
\proof
Applying Fact \ref{fact-app} (1), there is a definable open $V\ni a$ such that $a\in V\sub X_b$ and $\dim(b/A[V])=\dim(b/A)$.

Because $b$ is generic in $K^n$ over $A[V]$, it remains generic in $Y=\{b'\in K^n: V\sub X_{b'}\}$. it follows that $\dim(Y)=n$ and  $b\in Int(Y)$, so
$(a,b)\in V\times Int(Y)\sub Int(X).$\qed

To prove the lemma, let $V\sub K^r$ be an open neighborhood of $g(a,b)$. By Fact \ref{fact-app} (1), we may replace $V$ by $V_1$, $g(a,b)\ni V_1\sub V$, with $b$ generic in $K^n$ over
$A[V_1]$.  Consider the set $X=\{(x,y)\in W: f(x,y)\in V_1\}$. We need to see that $(a,b)\in Int(X)$. Since $f(x,b)$ is continuous at $a$, we have $a\in Int(X^b)$, and hence by the
above claim, $(a,b)\in Int(X)$. \qed

We now return to the proof of Proposition \ref{appendix-prop}.  Just like in  \cite{F-K}, we first reduce to the case where $a=0\in K^m$ and $g(0,y)\equiv 0$.

 After permuting $a$, we may write it as $a=(a_1,a_2)$ where $(a_2,b)$ is $acl$-independent over $A$ and $a_1\in acl(a_2bA)$. The set $W^{(a_2,b)}$ is open and $f(x_1,a_2,b)$ is still
 $C^1$ at $a_1$. Thus, by replacing $b$ with $(a_2,b)$ and $a$ with $a_1$, we may assume that $a\in acl(bA)$. By Fact \ref{appendix-cor}, we may add parameters to $A$ while preserving
 the genericity of $b$, such that $a\in dcl(bA)$. We still use $A$ for this new parameter set. Thus,  $b=\alpha(a)$ for an $A$-definable function $\alpha$. Since $b$ is generic in $dom(\alpha)$, then $\alpha$ is continuously
 differentiable at $b$. Without loss of generality, $dom(\alpha)=\pi_1(W)$.

Consider the local $C^1$-diffeomorphism $\bar \alpha:(x,y)\mapsto (x-\alpha(y),y)$. It sends $W$ to a set $\bar W$ and $(a,b)$ to $(0,b)$, so by Fact \ref{appendix-claim}, $(0,b)\in Int(\bar
W)$. The pushforward of $g$ via $\bar \alpha$ is $\bar g(x,y)=g(x+\alpha(y),y)$. The map $\bar g(x,b)$ is still $C^1$ on $\bar W^b$, and it is sufficient to prove that $\bar g$ is $C^1$ at
$(0,b)$. So, we may replace $g$ with $\bar g$, $W$ with $\bar W$ and $(a,b)$ with $(0,b)$. We still use $g$ and $W$ for the sets. Finally, since $b$ is generic in $K^n$, it follows from
Fact \ref{fact-app} (3) that  the function $g(0,y)$ is $C^1$ in a neighborhood of  $b$, so we may replace $g$ with $g(x,y)-g(0,y)$, thus assume that $g(0,y)\equiv 0$, and in particular,
$D_yg_{(0,b)}=0$, so $Dg_{(0,b)}=(D_xg_{(0,b)},0)\in M_{r\times (m+n)}(K)$.

To simplify notation below, we view $x\in K^n$ both as a row and a column vector, depending on context. Thus, for, say, $(x,y)\in K^m\times K^n$, we write  $Dg_{{(a,b)}}\cdot (x,y)$,
instead of $Dg_{(a,b)}(x,y)^t$.

Notice that in order to show that $g$ is diffenrentiable at $(0,b)$ we need to show that for every $\epsilon\in \Gamma$, the point $(0,b)$ belongs to the interior of the set of
$(x,y)\in K^m\times K^n$, such that
$$\|g(x,y)-g(0,b)-Dg_{(0,b)}\cdot (x-0,y-b)\|<\epsilon \|(x,y-b)\|,$$ which, since $g(0,b)=0$ and $D_yg_{(0,b)}=0$,  equals
\begin{equation}\label{appendix-eq1}\left\{(x,y)\in K^m\times K^n: \|g(x,y)-D_xg_{(0,b)}\cdot x\|<\epsilon \|(x,y-b)\|\right\}.\end{equation}

We fix $\epsilon\in \Gamma$.
By our assumption that $g(x,b)$ is differentiable at $0$, it follows that $0$ is in the interior of
$$\left\{x\in K^m: \|g(x,b)-D_xg_{(0,b)}\cdot x\|<\epsilon \|x\|\right\}.$$

By Claim \ref{appendix-claim}, $(0,b)$ is in the interior of
$$\left\{(x,y)\in K^m\times K^n:\|g(x,y)-D_xg_{(0,y)}\cdot x\|<\epsilon\|x\|\right\},$$ hence there exists $\delta_1\in \Gamma$ such that if $\|(x,y-b)\|<\delta_1$ then

$$\|g(x,y)-D_xg_{(0,y)}\cdot x\|<\epsilon\|x\|.$$

In order to prove that $(0,b)$ is in the interior of the set in (\ref{appendix-eq1}), we write \begin{equation}\label{appendix-eq2}g(x,y)-D_xg_{(0,b)}\cdot x=g(x,y)-D_xg_{(0,y)}\cdot
x+(D_xg_{(0,y)}-D_xg_{(0,b)})\cdot x.\end{equation}

\begin{claim} There is $\delta_2\in \Gamma$, such that for all $(x,y)\in K^m\times K^n$, if $\|y-b\|<\delta_2$ then $$\|(D_xg_{(0,y)}-D_xg_{(0,b)})\cdot x\|<\epsilon\|x\|.$$
\end{claim}
\proof We first observe that for every $A=(a_{i,j})\in M_{n}(K)$, if for all $i,j$, $|a_{i,j}|<\epsilon$ then for all $x\in K^n$, we have $\|A\cdot x\|<\epsilon\|x\|.$

Consider the map $G:K^n\to M_{r\times n}(K)$, $G(y)=D_xg_{(0,y)}$. It is definable over $A$ and hence continuous at $b$.
Thus, there exists $\delta_2\in \Gamma$ such that whenever $\|y-b\|<\delta_2$, then $\|D_xg_{(0,y)}-D_xg_{(0,b)}\|<\epsilon.$ The result follows from our above observation.\qed

If we now take $\delta=\min\{\delta_1,\delta_2\}$, for $\delta_2$ as in the above claim, then for all $(x,y)\in K^m\times K^n$ with $\|(x,y-b)\|<\delta$, we have, using (\ref{appendix-eq2}),

$$\|g(x,y)-D_xg_{(0,b)}\cdot x\|\leq
\max\{\|g(x,y)-D_xg_{(0,y)}\cdot x\|,\|(D_xg_{(0,y)}-D_xg_{(0,b)})\cdot x\|\}$$ $$<\epsilon\|x\|\leq \epsilon\|(x,y-b)\|
.$$

This ends the proof that $g(x,y)$ is differentiable at $(a,b)$. Since $a\in W^b$ was arbitrary it follows that for all $x\in W^b$, $g(x,y)$ is differentiable at $(x,b)$.
Consider now the map $G:(x,y)\mapsto Dg_{(x,y)}$. Since $g(x,b)$ is $C^1$ on $W^b$, the map $G(x,b)$ is continuous on $W^b$, and therefore by Lemma \ref{appendix-lemma1},  $G$ is continuous at $(a,b)$. Thus, $g$ is $C^1$ at
$(a,b)$. This ends the proof of Proposition \ref{appendix-prop}.
\qed


\end{document}